\documentclass[12pt]{article}
\usepackage{amsthm,amssymb,amsmath}
\usepackage[american,british]{babel}

\newtheorem*{Th}{Theorem}

\newcommand{\bX}{\boldsymbol X}

\newcommand{\bx}{\boldsymbol x}

\newcommand{\bb}{{\boldsymbol b}}

\newcommand{\bv}{{\boldsymbol v}}
\newcommand{\bomega}{{\boldsymbol \xi}}

\newcommand{\R}{{\mathbb R}}

\newcommand{\D}{{\partial}}

\begin{document}

\title{Sequences of Levy Transformations and \\
Multi-Wro\'nski Determinant Solutions \\of  the Darboux System}

\author{Q. P. Liu\thanks{On leave of absence from
Beijing Graduate School, CUMT, Beijing 100083, China}
\thanks{Supported by {\em Beca para estancias temporales
de doctores y tecn\'ologos extranjeros en
Espa\~na: SB95-A01722297}}
   $\,$ and Manuel Ma\~nas\thanks{Partially supported by CICYT:
 proyecto PB95--0401}\\
 Departamento de F\'\i sica Te\'orica,
\\ Universidad Complutense,\\
E28040-Madrid, Spain.}

\date{}

\maketitle

\begin{abstract}
Sequences of Levy transformations for the Darboux system of
conjugates nets in multidimensions are studied. We show that after a suitable 
number of Levy
transformations, with at least a Levy transformation in each direction,
we get closed formulae in terms of multi-Wro\'nski determinants.
These formulae are for the tangent vectors, Lam\`e coefficients,
rotation coefficients and points of the surface.
\end{abstract}
\newpage


{\bf 1}. The interaction between Soliton Theory and Geometry is a growing
subject. In fact, many systems that appear by geometrical
considerations have been studied independently in Soliton Theory,
well-known examples include the Liouville and sine-Gordon equations which
characterize minimal and pseudo-spherical surfaces, respectively.
Another relevant case is given by the the Darboux equations that
were solved 12 years ago in its matrix generalization, using the
$\bar\partial$--dressing, by Zakharov and Manakov \cite{zm}.

In this note we want to iterate a transformation that preserves the
Darboux equations which is known as Levy transformation
\cite{levy}.

{\bf 2}. The Darboux equations
\begin{equation} \label{dar}
\frac{\partial\beta_{ij}}{\D u_k}=\beta_{ik}\beta_{kj},\;\;
 i,j,k=1,\dotsc, N,\quad i\ne j\ne k\ne i,
\end{equation}
for the $N(N-1)$ functions
$\{\beta_{ij}\}_{\substack{i,j=1,\dotsc,N\\i\neq j}}$ of
$u_1,\dotsc,u_N$, characterize $N$-dimen\-sional submanifolds of
$\R^P$, $N\leq P$, parametrized by conjugate coordinate systems
\cite{Darboux,Eisenhart}, and are the compatibility conditions of
the following linear system
\begin{equation} \label{X}
\frac{\partial \bX_j}{\D u_i} = \beta_{ji} \bX_i, \quad i,j=1,\dotsc,N,\quad
i\ne j,
\end{equation}
involving suitable $P$-dimensional vectors $\bX_i$,
tangent to the coordinate lines. The so called Lam\'e coefficients
satisfy
\begin{equation} \label{H}
\frac{\partial H_j}{\D u_i} = \beta_{ij} H_i, \quad i,j=1,\dotsc,N,\quad
i\ne j,
\end{equation}
and the points of the surface $\bx$ can be found by means of
\begin{equation}\label{points}
\frac{\D \bx}{\D u_i}= \bX_i H_i,\quad i=1,\dotsc, N,
\end{equation}
which is equivalent to the more standard Laplace equation
\[
\frac{\D^2\bx}{\D u_i\D u_j}=\frac{\D \ln H_i}{\D u_j}\frac{\D \bx}{\D u_i}
+\frac{\D \ln H_j}{\D u_i}\frac{\D \bx}{\D u_j},
\quad i,j=1,\dotsc,N,\;\; i\neq j.
\]

A Darboux type transformation for this system was found by
Levy \cite{levy,eisenhart,ks}. In  fact, in \cite{levy} the
transformation is constructed only for two-dimensional surfaces,
$N=2$, being the Darboux equations in this case trivial and 
Levy only presents the transformation for the points of the surface.
However, in \cite{ks} the Levy transformation is extended to the
first non trivial case of Darboux equations, namely $N=3$.
The extension to arbritary $N$ is straightforward and reads as follows.
Given a solution $\xi_j$ of
\[
\frac{\D \xi_j}{\D u_k}=\beta_{jk} \xi_k,
\]
for each of the
$N$ possible directions in the
coordinate space there is a corresponding Levy transformation
that reads for the $i$-th case:
\begin{align*}
&\bx[1]=\bx-\frac{\Omega(\xi,H)}{\xi_i }{\bX_i},
\\[4mm]
&\begin{cases}
\bX_i[1]=\dfrac{1}{\xi_i}\Big(\xi_i\dfrac{\D\bX_i}{\D u_i}-
\dfrac{\D\xi_i}{\D u_i}\bX_i\Big),\\[3.5mm]
\bX_k[1]=\dfrac{1}{\xi_i}(\xi_i\bX_k-\xi_k\bX_i),
\end{cases}
\\[4mm]
&\begin{cases}
H_i[1]=-\dfrac{\Omega(\xi,H)}{\xi_i},\\
H_k[1]=H_k-\beta_{ik}\dfrac{\Omega(\xi,H)}{\xi_i},
\end{cases}
\\[4mm]
&
\begin{cases}
\beta_{ik}[1]=-\dfrac{1}{\xi_i}\Big(\beta_{ik}\dfrac{\D \xi_i}{\D u_i}-
\xi_i\dfrac{\D\beta_{ik}}{\D u_i}\Big),\\[3.5mm]
\beta_{ki}[1]=-\dfrac{\xi_k}{\xi_i},\\[3.5mm]
\beta_{kl}[1]=-\dfrac{\xi_k\beta_{il}-\xi_i\beta_{kl}}{\xi_i},
\end{cases}
\end{align*}
where $k,l=1,\dotsc, N$ with $k\neq l \neq i$. Here we have
introduced the potential $\Omega(\xi,H)$ defined by
\[
\frac{\D\Omega(\xi,H)}{\D u_k}=\xi_kH_k,\quad k=1,\dotsc,N,
\]
which are compatible equations by means of the equations
satisfied by $\xi_k$ and $H_k$.

{\bf 3}. Using  Crum type ideas \cite{crum} one can iterate this
Levy transformation. However now there is a difference with
respect to the iteration of the Darboux transformation of the
1-dimensional Schr\"odinger equation: we have $N$ different
elementary Levy's transformations $\{{\cal L}_i\}_{i=1,\dotsc,N}$.

If one performs  less than $N$ iterations or more than $N$
iterations, say ${\cal L}_{i_1}\dotsb{\cal L}_{i_M}$ with
$\{1,\dotsc,N\}\not\subset\{i_1,\dotsc,i_M\}$, one gets non
symmetric formulae in which the initial $\beta$'s and its
derivatives appear explicitly. However, if in the latter case we
have $\{1,\dotsc,N\}\subset\{i_1,\dotsc,i_M\}$, that is we have
perform at least one Levy transformation in each spatial direction
we obtain formulae only in terms of Wro\'nski determinants of the
wave functions with no $\beta$'s appearing explicitly.

To present our main result, we introduce some convenient notations.
First we define $\D_i:=\D/\D u_i$. Second, for any set of functions
 $\{\xi^i_j\}_ {\substack{i=1,\dots,M\\ j=1.\dots, N}}$ we denote by
$W_j(n)$ the following Wro\'nski matrix
\[
W_j(n):= W_j(\xi^1_j, \dotsc, \xi^M_j):=\begin{pmatrix}
\xi^1_j&\xi^2_j&\dots&\xi^M_j\\
\D_j\xi^1_{j}&\D_j\xi^2_{j}&\dots&\D_j\xi^M_{j}\\
\vdots&\vdots& &\vdots\\
\D^{n-1}_j\xi^1_{j}&\D^{n-1}\xi^2_{j}&\dots&\D^{n-1}_j\xi^M_j\\
\end{pmatrix}.
\]
For any partition of $M=m_1+m_2+ \dots + m_N$, we construct
a multi-Wro\'nski matrix
\[
{\cal W}:=\begin{pmatrix}W_1(m_1)\\ W_2(m_2)\\ \vdots \\
W_N(m_N)\end{pmatrix}.
\]

Now we are ready to present the following:
\begin{Th}\label{theorem}
Given $M$ functions $\{\xi_i^j\}_{\substack{i=1,\dotsc,N\\ j=1,\dotsc,
 M}}$ and
$\bX_i=(X_i^1,\dotsc,X_i^P)^{\operatorname{t}}$, $i=1,\dotsc,N$, all of them
solutions of \eqref{X} and $H_i$, $i=1,\dotsc,N$, solutions
of \eqref{H},
for  given $\beta_{ij}$,  then new solutions
$\bX_i[M],H_i[M]$ and $\beta_{ij}[M]$ are defined by:
\[
X_i^\ell[M]=\frac{\left|{\mathbb X}_i^\ell\right|}{\left|{\cal
W}\right|},\quad H_i[M]=-\frac{\left|{\mathbb
H}_i\right|}{\left|{\cal W}\right|},\quad
\beta_{ij}[M]=-\frac{\left|{\cal W}_{ij}\right|}{\left|{\cal W}\right|},
\]
where
\[
{\mathbb X}_i^\ell=\begin{pmatrix}{\cal W}&\bv^\ell\\
\D_i^{m_i}\bomega_i&\D_i^{m_i}X_i^\ell\end{pmatrix},
\]
with
\begin{align*}
\bv^\ell&:=(\bv_1^\ell,\dotsc,\bv_N^\ell)^{\operatorname{t}},\text{ being }
\bv_k^\ell:=(X_k^\ell,\D_k X_k^\ell,\dotsc,
\D_k^{m_k-1}X_k^\ell),\\
\bomega_i&:=(\xi_i^1,\dotsc,\xi_i^M),
\end{align*}
${\mathbb H}_i$ is obtained from $\cal W$ by replacing the last row
of the $i$-th block by $\Omega(\bomega,H)$ and ${\cal W}_{ij}$
by replacing the last row of the $j$-th block by  $\D_i^{m_i}\bomega_i$.
In the partition $M=m_1+m_2+ \dots + m_N$ we need $m_i\in\mathbb N$.

Moreover, for the transformed surface we have the parametrization
\[
\bx[M]=\frac{1}{|{\cal W}|}\bigg(\begin{vmatrix}{\cal W} &\bv^1\\
\Omega(\bomega,H) &x^1
\end{vmatrix},\dotsc,\begin{vmatrix}{\cal W} &\bv^P\\
\Omega(\bomega,H) &x^P
\end{vmatrix} \bigg)^{\operatorname{t}}.
\]
\end{Th}

\begin{proof}
 The proof that follows is inspired by  \cite{freeman,nimmo},
however is extended to this multicomponent system and we
give a more detailed account.

We first need to show that
\[
\D_kX_i^\ell[M]=\beta_{ik}[M]X_k^\ell[M],
\]
or equivalently that the following bilinear equation holds
\[
\left|{\cal W}\right|\D_k\left|{\mathbb X}_i^\ell\right|-
\left|{\mathbb X}_i^\ell\right| \D_k\left|{\cal W}\right|+
\left|{\mathbb X}_k^\ell\right|\left|{\cal W}_{ik}\right|=0.
\]
To this aim we consider the following $(2M+1)\times(2M+1)$ square
matrix
\[
{\cal A}_{ik}^\ell:=
\begin{pmatrix}
\quad A_k\quad& 0 &\D_k^{m_k-1}\bomega_k^{\operatorname{t}}&
\D_k^{m_k}\bomega_k^{\operatorname{t}}&
\D_i^{m_i}\bomega_i^{\operatorname{t}}\\[3mm]
\quad 0\quad &\quad A_k\quad&\D_k^{m_k-1}\bomega_k^{\operatorname{t}}&
\D_k^{m_k}\bomega_k^{\operatorname{t}}&
\D_i^{m_i}\bomega_i^{\operatorname{t}}\\[3mm]
0&\bb_k^\ell&\D_k^{m_k-1}X_k^\ell&\D_k^{m_k}X_k^\ell&\D_i^{m_i}X_i^\ell
\end{pmatrix},
\]
where $A_k$ is a $M\times(M-1)$ rectangular matrix
\[
(A_k)^{\operatorname{t}}:=
\begin{pmatrix}
 W_1(m_1)\\
\vdots\\
\hat W_k(m_k)
\\
\vdots\\
W_N(m_N)
\end{pmatrix},
\]
with $\hat W_k(m_k)$ obtained from $W_k(m_k)$ by deleting the last row,
and
\[
\bb_k^\ell=(\bv_1^\ell,\cdots,\hat\bv_k^\ell,\dotsc,\bv_N^\ell),
\]
with $\hat\bv_k^\ell$ obtained by deleting the last element in $\bv_k^\ell$.

We now recall the Laplace's general expansion theorem \cite{algebra}
that we shall use in this proof, this theorem allows the computation
of an $n\times n$ matrix $A:=(a_{ij})$ as follows:
\begin{multline*}
\det A=\sum_{\substack{\rho_1,\dotsc,\rho_r\\
\rho_1<\dotsb<\rho_r}}(-1)^{\gamma_1+\dotsb+\gamma_r+\rho_1+\dotsb+\rho_r}\\
\begin{vmatrix}
a_{\gamma_1\rho_1}&a_{\gamma_1\rho_2}&\hdots &a_{\gamma_1\rho_r}\\
a_{\gamma_2\rho_1}&a_{\gamma_2\rho_2}&\hdots &a_{\gamma_2\rho_r}\\
\vdots &\vdots&&\vdots\\
a_{\gamma_r\rho_1}&a_{\gamma_r\rho_2}&\hdots &a_{\gamma_r\rho_r}
\end{vmatrix}\times
\begin{vmatrix}
a_{\delta_1\sigma_1}&a_{\delta_1\sigma_2}&\hdots &a_{\delta_1\sigma_s}\\
a_{\delta_2\sigma_1}&a_{\delta_2\sigma_2}&\hdots &a_{\delta_2\sigma_s}\\
\vdots &\vdots&&\vdots\\
a_{\delta_s\sigma_1}&a_{\delta_s\sigma_2}&\hdots &a_{\delta_s\sigma_s}
\end{vmatrix},
\end{multline*}
where $r+s=n$ and
\begin{align*}
(\gamma_1,\dotsc,\gamma_r,\delta_1,\dotsc,\delta_s)&=(1,\dots,n)\\
(\rho_1,\dotsc,\rho_r,\sigma_1,\dotsc,\sigma_s)&=(1,\dots,n),
\end{align*}
up to permutations.

Let us now expand the determinant of the matrix ${\cal A}^\ell_{ik}$ by means of
the Laplace's general expansion theorem. Here, we take
$r=M$, $\gamma_i=i$ and
$\delta_i=M+i(i=1,\dotsc,M)$. It is easy
to see that:
\begin{align*}
\left| {\cal A}_{ik}^\ell\right|=(-1)^{M-1}&\Bigg(
\begin{vmatrix}A_k&\D^{m_k-1}_{k}\bomega_k^{\operatorname{t}}
\end{vmatrix}\times
\begin{vmatrix}A_k&\D_k^{m_k}\bomega_k^{\operatorname{t}}&
\D_i^{m_i}\bomega_i^{\operatorname{t}}\\
\bb_k^\ell&\D_k^{m_k}X_k^\ell&\D_i^{m_i}\bomega_i
\end{vmatrix}\\
&-
\begin{vmatrix}A_k&\D^{m_k}_{k}\bomega_k^{\operatorname{t}}
\end{vmatrix}\times
\begin{vmatrix}A_k&\D_k^{m_k-1}\bomega_k^{\operatorname{t}}&
\D_i^{m_i}\bomega_i^{\operatorname{t}}\\
\bb_k^\ell&\D_k^{m_k-1}X_k^\ell&\D_i^{m_i}\bomega_i
\end{vmatrix}\\
&+
\begin{vmatrix}A_k&\D^{m_i}_{i}\bomega_i^{\operatorname{t}}
\end{vmatrix}\times
\begin{vmatrix}A_k&\D_k^{m_k-1}\bomega_k^{\operatorname{t}}&
\D_i^{m_k}\bomega_i^{\operatorname{t}}\\
\bb_k^\ell&\D_k^{m_k-1}X_k^\ell&\D_k^{m_k}\bomega_k
\end{vmatrix}
\Bigg),
\end{align*}
expression that after an even number of permutations of columns
and transposition reads
\[
\big|{\cal A}_{ik}^\ell\big|=
(-1)^{M-1}\big[|{\cal W}|\D_k|{\mathbb X}_i^\ell|-|{\mathbb
X}_i^\ell|
\D_k|{\cal W}|
+|{\mathbb X}_k^\ell||{\cal W}_{ik}|\big].
\]
But the Laplace's theorem also implies $|{\cal A}^\ell_{ik}|=0$, to see this
we just use the standard version of this theorem and expand the determinant
with respect to its last row. In doing so we get a sum in which all terms
vanish, this last statement follows again from the Laplace's general
expansion theorem.
This gives the desired result.

Next we prove that
\[
\D_kH_i[M]=\beta_{ki}[M]H_k[M],
\]
or equivalently that the following bilinear equation holds:
\[
{\left|\cal W\right|}\D_k\left|{\mathbb H}_i\right|-
\left|{\mathbb H}_i\right|\D_k\left|{\cal W}\right|
+\left|{\cal W}_{ki}\right|\left|{\mathbb H}_k\right|=0.
\]
As before this relation is a consequence of the Laplace's general
expansion theorem. For this aim we consider the  $2M\times 2M$ square
matrix:
\[
{\cal B}_{ik}:=
\begin{pmatrix}
\quad B_{ik}\quad& 0  &
\D_k^{m_k-1}\bomega_k^{\operatorname{t}}&
\D_k^{m_k}\bomega_k^{\operatorname{t}}&
\D_i^{m_i-1}\bomega_i^{\operatorname{t}}&\Omega(\bomega,H)^{\operatorname{t}}\\[
3mm]
\quad 0\quad &\quad B_{ik}\quad
&\D_k^{m_k-1}\bomega_k^{\operatorname{t}}&
\D_k^{m_k}\bomega_k^{\operatorname{t}}&
\D_i^{m_i-1}\bomega_i^{\operatorname{t}}&\Omega(\bomega,H)^{\operatorname{t}}
\end{pmatrix},
\]
where $B_{ik}$ is a $M\times(M-2)$ rectangular matrix
\[
(B_{ik})^{\operatorname{t}}:=
\begin{pmatrix}
 W_1(m_1)\\
\vdots\\
\hat W_i(m_i)\\
\vdots\\
\hat W_k(m_k)
\\
\vdots\\
W_N(m_N)
\end{pmatrix}.
\]
Using the version of the Laplace expansion appearing in \cite{freeman}
(eq. (3.3)) we get the desired bilinear formula.

Finally, we prove the formula for $x^\ell=\Omega(X^\ell[M],H[M])$ (see 
\eqref{points}).
 This is achived by considering the following $(2M+1)\times(2M+1)$ square matrix
\[
{\cal C}_{ik}^\ell:=
\begin{pmatrix}
\quad A_k\quad& 0 &\D_i^{m_i-1}\bomega_i^{\operatorname{t}}&
\D_i^{m_i}\bomega_i^{\operatorname{t}}&
\Omega(\bomega,H)^{\operatorname{t}}\\[3mm]
\quad 0\quad &\quad A_k\quad&\D_i^{m_i-1}\bomega_i^{\operatorname{t}}&
\D_i^{m_i}\bomega_i^{\operatorname{t}}&
\Omega(\bomega,H)^{\operatorname{t}}\\[3mm]
0&\bb_k^\ell&\D_i^{m_i-1}X_i^\ell&\D_i^{m_i}X_i^\ell&\Omega(X^\ell,H)
\end{pmatrix},
\]
and using  that $x^\ell=\Omega(X^\ell,H)$ and
Laplace's general expansion theorem.
\end{proof}

{\bf 4}. 
Sequences of Levy transformations for two dimensional surfaces
have already been studied in \cite{hammond,t}, see also \cite{eisenhart}. 
Let us remark that the  Darboux equations are trivial in this case
and that they only consider the points in the surface.
Up to a factor  ($(H_1\dotsb H_N|{\cal W}|)^{-1}$) and the
choice $H_i\xi_i^j=\D_i\theta^{(j)}$ our formula
for the points of the surface coincides, when $N=2$, with
the formula of \cite{hammond}, in where, to our knowldege,
is the first place where double wronskian appeared.

From a complete different point of view Nimmo considered in \cite{nimmo}
what he called Darboux transformations for the two-dimensional
Zakharov-Shabat/ AKNS spectral problem, i. e.
in the context of the Davey-Stewartson equations. 
In fact, this is intimately 
connected with two-dimensional conjugate nets \cite{k}.
His results are special cases of ours: first we have arbitary dimension,
not only $N=2$ as in \cite{nimmo}; second our partition for $N=2$,
$M=m_1+m_2$,
is more general than his, $M=2m$; third
we have computed not only the
transformation for the potentials $\beta_{12}=q$ and $\beta_{21}=r$ and
wave functions but also for the adjoint wave function and for the 
corresponding points in the surface, that in this case, as we mentioned in
the previous paragraph can be found in \cite{hammond}.

The above remarks illustrate the fact that same problem has been
tackled by different techniques coming form Geometry on one hand
and Soliton Theory on the other, covering different aspects of it.
In this paper we have extended the results of both approaches
to higher dimensions,
in where the Darboux equations are not any more trivial.
In fact, from the Soliton Theory point of view the Levy transformation
for the Darboux system can be consiedered as a {\em elementary Darboux
transformations} for the $N$-component Kadomtsev-Petviashvili hierarchy
\cite{djkm}.

We already mentioned that there exist other possibilities for
iterations. In fact we have requested that there is at least one
 Levy transformation per
direction. If this is not the case our results do not hold any more.
However, one could get closed formulae in where the $\beta$'s
appear explicitly. In principle, one has $N$ possible {\em different}
types of formulae. But we are not going to consider
this problem in this Letter.


\begin{thebibliography}{99}

\bibitem{crum} M. Crum, {\em Quart. J. Math.} {\bf 6} (1955) 121.

\bibitem{djkm} E. Date, M. Jimbo, M. Kashiwara, and T. Miwa,
{\em J. Phys. Soc. Japan} {\bf 50} (1981) 3806.

\bibitem{Darboux}
G. Darboux,
{\it Le\c{c}ons sur la th\'{e}orie g\'{e}n\'{e}rale des surfaces IV, Liv. VIII,
Chap. XII}, Chelsea Publishing Company,
New York (1972).

\bibitem{Eisenhart}
L. P. Eisenhart,
{\it A Treatise on the Differential Geometry of Curves and Surfces},
Ginn and Co., Boston (1909).

\bibitem{eisenhart} L. P. Eisenhart, {\em Transformations
of Surfaces}, Chelsea Publishing Company, New York (1962).

\bibitem{freeman} N.C. Freeman, {\em IMA J. Appl. Math.} {\bf 32}
(1984) 125.


\bibitem{hammond} E. S. Hammond, {\em Ann. Math.} {\bf 22} (1920) 238.

\bibitem{k} B. G. Konopelchenko, {\em Phys. Lett.} {\bf A183} (1993) 153.

\bibitem{ks} B. G. Konopelchenko and W. K. Shief,
{\em Lam\'e and Zakharov-Manakov systems: Combescure, Darboux and
B\"acklund transformations}, Preprint AM93/9, UNSW (1993).

\bibitem{levy} L. Levy, {\em J. l'\'Ecole Polytecnique} {\bf 56} (1886) 63.

\bibitem{nimmo} J. J. C. Nimmo, {\em Inverse Problems} {\bf 8} (1992) 219.

\bibitem{algebra} O. Schreier and E. Sperner, {\em Introduction to
Modern Algebra and Matrix Theory}, Chelsea Publishing Company,
New York (1951).

\bibitem{t} G. Tzitzeica, {\em C. R. Acad. Sci. Paris} {\bf 156} (1913) 375.

\bibitem{zm}
V. E. Zakharov and S. E. Manakov,
{\em Func. Anal. Appl.} {\bf 19} (1985) 11.
\end{thebibliography}
\end{document}